\documentclass[a4paper]{article}
\usepackage{amsmath,amssymb} 
\usepackage{diags}
\usepackage[amsmath,thmmarks]{ntheorem}
\usepackage{bm} 
\usepackage{graphicx}
\usepackage{float}
\usepackage{overpic}
\usepackage{esvect}
\usepackage{epstopdf}
\usepackage{color}
\usepackage[urlcolor=blue]{hyperref}

\theoremstyle{plain}
\newtheorem{theorem}{Theorem}
\newtheorem{proposition}{Proposition}
\newtheorem*{thm}{Theorem}

\theoremstyle{definition}
\theorembodyfont{\normalfont}
\newtheorem{lemma}{Lemma}

\theoremstyle{remark}

\newarrow{DL}{=}{=}{=}{=}{=}
\title{Atoroidal Manifolds in Small Covers}
\author{Lisu Wu}
\date{}
\begin{document}
\maketitle
\begin{abstract}
We show that a $3$-dimensional small cover is atorodal if and only if there is no $4$-belt in the corresponding simple polytope.
\end{abstract}
\section{Introduction}
Small covers are a class of closed manifolds which admit locally standard $\mathbb{Z}_2^n$-actions, such that the orbit spaces are some simple convex ploytope in an Euclidean space. Or defined equivalently, an $n$ dimensional small cover consists of $2^n$ copies of a simple convex ploytope whose faces are guled according to the colorings in $\mathbb{Z}_2^n$. In general,
the topological properties of a small cover are closely related to the combinatorics of the orbit ploytope and the colorings on its facets(codimension-$1$ faces). For instance, the cohomlogy ring of a small cover with $\mathbb{Z}_2$ coefficients is isomorphic to $\mathbb{Z}_2/(I+J)$, where $\mathbb{Z}_2/I$ is the Stanley-Reisner ring of the orbit ploytope, and $J$ depends on the colorings on its facets\cite[Theorem 4.14]{DJ1};  and a small cover is orientable if and only if the sum of entries of the coloring on each face  is odd\cite[Theorem 1.7]{NN}.

Let $P$ be an $n$-dimensional simple convex polytope in $\mathbb{R}^n$, and $\pi: M\longrightarrow P$ be a small cover over $P$.  We denote the set of all facets of $P$ by $\mathcal{F}(P)$, 
then there is a {\em characteristic map} $\lambda: \mathcal{F}(P)\rightarrow \mathbb{Z}_2^n$ with the coloring of each facet identical with $\lambda(F)$, where $\mathbb{Z}_2^n$ is taken as a multiplicative group. For any proper face $f$ of $P$, we define $$G_f\overset{\Delta}{=}\text{the subgroup of } \mathbb{Z}_2^n \text{ generated by the set} \{\lambda(F)\mid f\subset F\}.$$
and $G_P\overset{\Delta}{=}\{1\}\subset \mathbb{Z}_2^n$. Then $M$ is homeomorphic to the quotient space 
\begin{equation}
P\times \mathbb{Z}_2^n/\sim
\end{equation}
where $(p,g)\sim(q,h)$ if and only if $p=q$ and $g^{-1}h\in G_{f(p)}$, and $f(p)$ is the unique face of $P$ that contains $p$ in its relative interior.
The {\em right-angled Coxeter group} of $P$ is defined by:
\begin{equation}
W_P\overset{\Delta}{=}\langle s_F, F \in \mathcal{F}(P)\mid s_F^2=1, \forall F; (s_Fs_{F'})^2=1, \forall F\cap F'\neq\varnothing \rangle
\end{equation}
By giving the Borel construction $M_{\mathbb{Z}_2^n}=E\mathbb{Z}_2^n\times_{\mathbb{Z}_2^n}M$ of $M$, Davis and Januszkiewicz \cite[Corllory 4.5]{DJ1} proved that the fundamental group of a small cover, denoted as $\pi_1(M)$, is isomorphic to the kernel of a homomorphism from $W_P$ to $Z_2^n$, which is induced by the characteristic map over $\mathcal{F}(P)$.  Actually, there is a short right split group sequence
\begin{equation}
1\longrightarrow\pi_1(M)\longrightarrow W_P\overset{\phi}{\longrightarrow}\mathbb{Z}_2^n\longrightarrow 1
\end{equation}
where $\phi(s_F)=\lambda(F)$ for each facet $F$ of $P$.
Wu and Yu \cite{WY} described this relation more explicitly based on the presentation of fundamental groups they calculated. Furthermore, they showed that a facial submanifold of a small cover is $\pi_1$-injective if and only of there is no $3$-belt conclude the corresponding face\cite[Theorem 3.3]{WY}. 

Let $P$ be a $3$-dimensional simple polytope. A {\em$k$-circuit} in $P$ is a simple loop on the boundary of $P$ which intersects transversely the interior of exactly $k$ distinct edges, and a $k$-circuit $c$ is called {\em prismatic} if the endpoints of all edges which $c$ intersects are distinct.
A {\em $k$-belt} in $P$ is a set of $k$ distinct faces $F_1, \cdots, F_k$ of $P$ such that $F_i\cap F_{i+1}\neq \varnothing$ for $1\leq i\leq k-1$,  $F_k\cap F_{1}\neq \varnothing$, and any three face in the belt have no common intersection.
It is clear that each $k$-belt determines a prismatic $k$-circuit. Furthermore, a prismatic $3$-circuit can determines a $3$-belt; and if there is no prismatic $3$-circuit, then a prismatic $4$-circuit determines a $4$-belt.
It is clearly that the cross section surrounded by a $4$-circuit $c$ determined by a $4$-belt is a square, we denote this cross section by $F$. If $c$ is prismatic, then $\pi^{-1}(F)$ is torus or Klein bottle in $M$, denoted by $M_F$.
And a closed $3$-manifold is {\em atoroidal} if it contains no essential torus, otherwise it is {\em toroidal}.

In this paper, we consider the $\pi_1$-injectivity of the sectional submanifold $M_F$ which is determined by a $4$-belt in simple polytope $P$. According to the following diagram, 
\begin{equation}
\begin{diagram}
1&\rTo&\pi_1(M_F)&\rTo^{\psi_F}&W_F&\rTo^{\phi_F}&\mathbb{Z}_2^3&&\\
&&\dTo^{i_*}&&\dTo_{j_*}&&\dDL&&\\
1&\rTo&\pi_1(M)&\rTo^{\psi}&W_P&\rTo^{\phi}&\mathbb{Z}_2^3&\rTo&1\\
\end{diagram}
\end{equation}
where $i_*$ and $j_*$ are induced by inclusion map of sectional submanifold $M_F$, we show that each $4$-belt in $P$ gives a $\pi_1$-injective torus or Klein bottle for $M$, i.e. $M$ is toroidal in this  case. 
 Furthermore, we give a topological proof for Andreev' Theorem in the right-angled case. The main result of this paper is the following theorem.
\begin{thm}
Let $M$ be a $3$-dimensional small cover over a simple polytope $P$, then $M$ is  atoroidal if and only if there is no $4$-belt in $P$. 
\end{thm}

The paper is organized as follows. In section 2, we construct a  presentation of  the fundamental group of the sectional submanifold which is determined by a $4$-belt. Using this presentation, we give a group homomorphism from $\pi_1(M_F)$ to the Coxeter group $W_F$.  In section 3, we show that sectional submanifold $M_F$ is $\pi_1$-injective. In section 4, we prove Andreev's Theorem in the right-angled case by the result in section 3. 
\newpage
\section{Presentations of  $\pi_1(M_F)$}
Let $P$ be a $3$-dimensional simple polytope, and $\pi: M\longrightarrow P$ be a small cover over $P$. If  $c$ is a prismatic $4$-circuit in $P$ determined by a $4$-belt $F_1,F_2,F_3,F_4$, where $F_1\cap F_3=F_2\cap F_4=\varnothing$. Then $c$ intersects exactly $4$ edges and $4$ faces of $P$.  Thus $c$ encloses a square in $P$, denoted by $F$. And four faces bounded $F$ are  $F_1,F_2,F_3,F_4$ respectively, the edge of $F$ denoted by $f_i=F\cap F_i$, for $i=1,2,3,4$.
  Then the Coxeter group of $F$  is
 \begin{equation}
\begin{split}
W_F\overset{\Delta}{=}\langle s_1, s_2, s_3, s_4\mid &s_1^2=s_2^2=s_3^2=s_4^2=1;\\&(s_1s_2)^2=(s_2s_3)^2=(s_3s_4)^2=(s_4s_1)^2=1\rangle
\end{split}
\end{equation}

 There are, essentially, $5$ cases of colorings on $\mathcal{F}=\{f_1,f_2,f_3,f_4\}$ which is induced by characteristic map $\lambda: \mathcal{F}(P)\longrightarrow \mathbb{Z}_2^3$. We define the induced characteristic map of $F$:
\begin{equation}
\lambda_F: \mathcal{F}\longrightarrow \mathbb{Z}_2^3
\end{equation}
 where $\lambda_F(\mathcal{F})=\lambda_F(f_1,f_2,f_3,f_4)\in \{(e_1,e_2,e_1,e_2),  (e_1,e_2,e_1,e_1e_2), (e_1,e_2,e_3,e_2),\\ (e_1,e_2,e_3,e_1e_2),  (e_1,e_2,e_3,e_1e_2e_3)\}$, and $e_1,e_2,e_3$ are basis of $\mathbb{Z}_2^3$.
 
Then the sectional submanifold  $M_F\overset{\Delta}{=}\pi^{-1}(F)$ determined by $F$ is
 \begin{equation}
M_F=F\times \mathbb{Z}_2^3/ \sim
\end{equation}
 where $(f_i, g)\sim (f_j, h)$ if and only if $i=j$ and $g^{-1}h=\langle \lambda_F(f_i)\rangle$. Thus sectional submanifolds determined by above $5$ colorings are double tori, double Klein bottles, torus, Klein bottle, torus, respectively.

 According to the construction of $M_F$, we consider the following two cases:
\\

{\bf Case 1:} $\lambda_F(\mathcal{F})=(e_1,e_2,e_1,e_2)  \text{ or } (e_1,e_2,e_1,e_1e_2)$.

In this case,  $\text{Im}(\lambda_F)$ generates the subgroup $\mathbb{Z}_2^2$ of $\mathbb{Z}_2^3$. The sectional submanifold is a disjoint union of two tori or two Klein bottles, one of which we denote by $M_F$. Then $M_F$ is glued by $4$ copies of $F$, we choose a vertex $p_0=f_1\cap f_2$ of $F$ as the base point of $\pi_1(M)$, and glue $\{(F,g)\}( g\in\mathbb{Z}_2^2)$ along its faces $\{(f_i,g)\}, i=1,2; g\in\mathbb{Z}_2^2$. For each face $(f_i, g)\subset (F,g)$, we choose a simple closed circles $\beta_{i,g}$, which crosses $p_0$ and $(f_i, g)$ in $M_F$,  as a generator of $\pi_1(M_F)$. Then 
 \begin{equation*}
 \begin{split}
 \pi_1(M_F, p_0)=\langle \beta_{i,g}, i=1,2,3,4; &\forall g\in \mathbb{Z}_2^2\mid \beta_{i,g}\beta_{i,g\cdot\lambda_F(f_i)}=1, \forall g;\\
 &\beta_{i,g}\beta_{j,g\cdot\lambda_F(f_i)}=\beta_{j,g}\beta_{i,g\cdot\lambda_F(f_j)}, \forall f_i\cap f_j\neq \varnothing, \forall g;\\
 &\beta_{i,g}=1, i=1,2; \forall g\rangle
  \end{split}
 \end{equation*}
 And consider the following short group sequence
 
  \begin{equation}\label{ES1}\begin{diagram}
1&\rTo&\pi_1(M_F,p_0)&\rTo^{\psi_F}&W_F& \pile{\rTo^{\phi_F} \\ \lTo_{\gamma_F}}&\mathbb{Z}_2^2&\rTo&1\\
\end{diagram} \end{equation}
where $\gamma_F: \mathbb{Z}_2^2\longrightarrow W_F$ is defined by $\gamma_F(e_1)=s_1, \gamma_F(e_2)=s_2$.
 And define
\begin{align*}
\psi_F:\pi_1(M_F,p_0)\longrightarrow& ~~W \\
 \beta_{i,g}\longmapsto& \gamma_F(g\cdot \lambda_F(f_i))s_i\gamma(g)\overset{\Delta}{=}S_{i,g}
\end{align*}
It is easy to check that both $\gamma_F$ and $\psi_F$ are well-defined. We have the following lemma.


\begin{lemma}
The sequence (\ref{ES1}) is right split and exact.
\end{lemma}
\noindent{\bf Proof.} cf. \cite[Lemma 2.9]{WY}  $\hfill{} \Box$
 \\
 
{\bf Case 2:}
 $\lambda_F(\mathcal{F})\in \{(e_1,e_2,e_3,e_2), (e_1,e_2,e_3,e_1e_2),  (e_1,e_2,e_3,e_1e_2e_3)\}$. 
 
 In this case,  $\text{Im}(\lambda_F)$ generates $\mathbb{Z}_2^3$. The sectional submanifold $M_F$ is glued by $8$ copies of $F$. Similarly, 
 we choose a vertex $p_0=f_1\cap f_2$ of $F$, and glue $\{(F,g)\}( g\in\mathbb{Z}_2^3)$ along along its faces $\{(f_i,g)\}$ for $i=1,2,  g\in\langle e_1,e_2\rangle$ or $i=3,  g\notin\langle e_1\rangle$. We shrink the faces $\{f_2, g\}, g\notin\langle e_1\rangle$ to a point, which is also denoted as $p_0$ and taken as the base point of $\pi_1(M)$. And for each pair $(f_i, g)\subset (F,g)$, we choose a simple closed circles $\beta_{i,g}$, which crosses $p_0$ and $(f_i, g)$ in $M_F$,  as a generator of $\pi_1(M_F)$. Then 
 \begin{equation}
 \begin{split}
 \pi_1(M_F)=\langle \beta_{i,g}, i=1,2,3,4; &\forall g\in \mathbb{Z}_2^3\mid \beta_{i,g}\beta_{i,g\cdot\lambda_F(f_i)}=1, \forall g;\\
 &\beta_{i,g}\beta_{j,g\cdot\lambda_F(f_i)}=\beta_{j,g}\beta_{i,g\cdot\lambda_F(f_j)}, \forall f_i\cap f_j\neq \varnothing; \forall g\\
 &\beta_{i,g}=1, i=1,2; \forall g \text{ or } i=3, \forall g \in G\rangle
  \end{split}
 \end{equation}
 where $G=\{e_1,e_1e_2,e_1e_3,e_1e_2e_3\}$.
 And consider the following short group sequence
  \begin{equation}\label{ES2}\begin{diagram}
1&\rTo&\pi_1(M_F,p_0)&\rTo^{\psi_F}&W_F& \pile{\rTo^{\phi_F} \\ \lTo_{\gamma_F}}&\mathbb{Z}_2^3&\rTo&1\\
\end{diagram} \end{equation}
where $\gamma_F: \mathbb{Z}_2^3\longrightarrow W_F$ is defined by $\gamma_F(e_1)=s_1,\gamma_F(e_2)=s_2,\gamma_F(e_3)=s_1s_3s_1$. And 
\begin{align*}
\psi_F:\pi_1(M_F,p_0)\longrightarrow& ~~W \\
 \beta_{i,g}\longmapsto& \gamma_F(g\cdot \lambda_F(f_i))s_i\gamma(g)\overset{\Delta}{=}S_{i,g}
\end{align*}
And both $\gamma_F$ and $\psi_F$ are well-defined. Then the similar lemma follows as above.
\begin{lemma}
The sequence (\ref{ES2}) is right split and exact. 
\end{lemma}
\noindent{\bf Proof.} $\phi_F\circ\gamma_F=id_{\mathbb{Z}_2^3}$ implies the right splitting of (\ref{ES2}).  We just need to prove $\pi_1(M_F,p_0)\cong \ker(\phi_F)$. $\phi_F\circ \psi_F(\beta_{i,g})=\phi_F(\gamma_F(g)s_i\gamma(g\cdot \lambda_F(f_i)))=\phi_F\circ\gamma_F(g)\cdot \phi_F(s_i)\cdot \phi_F\circ\gamma(g \lambda_F(f_i))=g\cdot\lambda_F(f_i)\cdot g \lambda_F(f_i)=1$, thus $\text{im}(\psi_F)\subseteq \ker(\phi_F)$.

Since $\text{im}(\psi_F)$ is a normal group of $W_F$, and
 \begin{equation}
\begin{split}
W_F&=\langle s_1, s_2, s_3, s_4\mid s_1^2=s_2^2=s_3^2=s_4^2=1;\\
&~~~~~~~~~~~~~~~~~~~~~~~(s_1s_2)^2=(s_2s_3)^2=(s_3s_4)^2=(s_4s_1)^2=1\rangle\\
&=\langle s_i, S_{i,g}, i=1,2,3,4. \forall g\in\mathbb{Z}_2^3\mid S_{i,g}S_{i,g\lambda_F(f_i)}=1,\forall i, g;\\
&~~~~~~~~~~~~~~~~~~~~~~~S_{i,g}S_{j,g\lambda_F(f_i)}=S_{j,g}S_{i,g\lambda_F(f_j)}, \forall f_i\cap f_j\neq \varnothing; \forall g\rangle\\
\end{split}
\end{equation}
thus $W_F/\text{im}(\psi_F)\cong \mathbb{Z}_2^3$ implies that $\text{im}(\psi_F)\cong\ker(\phi_F)$. Hence the sequence (\ref{ES2}) is exact.
 $\hfill{} \Box$
 
\section{$\pi_1$-injectivity of $M_F$}
In this section, we take $v=F_0\cap F_1\cap F_2$ as the base point, then the fundamental group of $M$ \cite[Proposition 2.1]{WY}
 \begin{equation}
 \begin{split}
 \pi_1(M,v)=\langle \alpha_{F,g}, \forall F\in\mathcal{F}(P); &\forall g\in \mathbb{Z}_2^3\mid \alpha_{F,g}\alpha_{F,g\cdot\lambda(F)}=1, \forall g;\\
 &\alpha_{F,g}\alpha_{F',g\cdot\lambda(F)}=\alpha_{F',g}\alpha_{F,g\cdot\lambda(F')}, \forall F\cap F'\neq \varnothing; \forall g\\
 &\alpha_{F,g}=1,\forall v\in F; \forall g \rangle
  \end{split}
 \end{equation}
and the Coxeter group of $P$
 \begin{equation}
 W_P=\langle t_F, F\in\mathcal{F}(P)\mid t_F^2=1, \forall F;(t_Ft_{F'})^2=1, for~ F\cap F'\neq \varnothing\rangle
 \end{equation}
Consider the following diagram
 \begin{equation}\label{Diag1}
\begin{diagram}
\pi_1(M_F)&\rTo^{\psi_F}&W_F\\
\dTo^{i_*}&&\dTo_{j_*}\\
\pi_1(M)&\rTo^{\psi}&W_P\\
\end{diagram}
\end{equation}
where $\psi$ is defined at \cite[Lemma 2.9]{WY}, and $s_0$ represents the generator in $W_P$ determined by the face $F_0$, $W_P\longrightarrow W_P/\langle s_0\rangle$ is a quotient map. $i_*$ and $j_*$ are induced by inclusion map.
\begin{align*}
j_*:W_F\longrightarrow&W_P \\
 s_i\longmapsto& t_{F_i}
\end{align*}
and in {\bf Case 1}, we define
\begin{align*}
i_*:\pi_1(M_F,p_0)\longrightarrow& \pi_1(M,v)\\
 \beta_{i,g}\longmapsto& \alpha_{F_i,g},~~g\in\langle e_1,e_2\rangle
\end{align*}
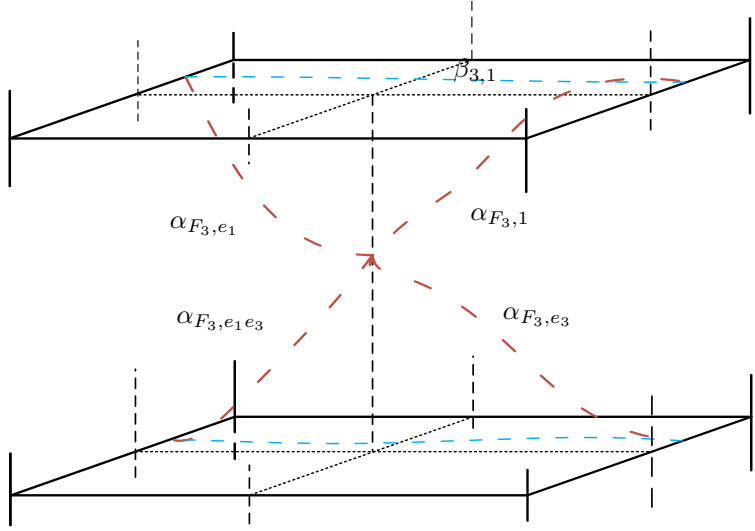
\begin{figure}[h]
\centering
\def\svgwidth{1.3\textwidth}
\begingroup%
  \makeatletter%
  \providecommand\color[2][]{%
    \errmessage{(Inkscape) Color is used for the text in Inkscape, but the package 'color.sty' is not loaded}%
    \renewcommand\color[2][]{}%
  }%
  \providecommand\transparent[1]{%
    \errmessage{(Inkscape) Transparency is used (non-zero) for the text in Inkscape, but the package 'transparent.sty' is not loaded}%
    \renewcommand\transparent[1]{}%
  }%
  \providecommand\rotatebox[2]{#2}%
  \newcommand*\fsize{\dimexpr\f@size pt\relax}%
  \newcommand*\lineheight[1]{\fontsize{\fsize}{#1\fsize}\selectfont}%
  \ifx\svgwidth\undefined%
    \setlength{\unitlength}{494.25bp}%
    \ifx\svgscale\undefined%
      \relax%
    \else%
      \setlength{\unitlength}{\unitlength * \real{\svgscale}}%
    \fi%
  \else%
    \setlength{\unitlength}{\svgwidth}%
  \fi%
  \global\let\svgwidth\undefined%
  \global\let\svgscale\undefined%
  \makeatother%
  \begin{picture}(1,0.4461305)%
    \lineheight{1}%
    \setlength\tabcolsep{0pt}%
    \put(0,0){\includegraphics[width=\unitlength,page=1]{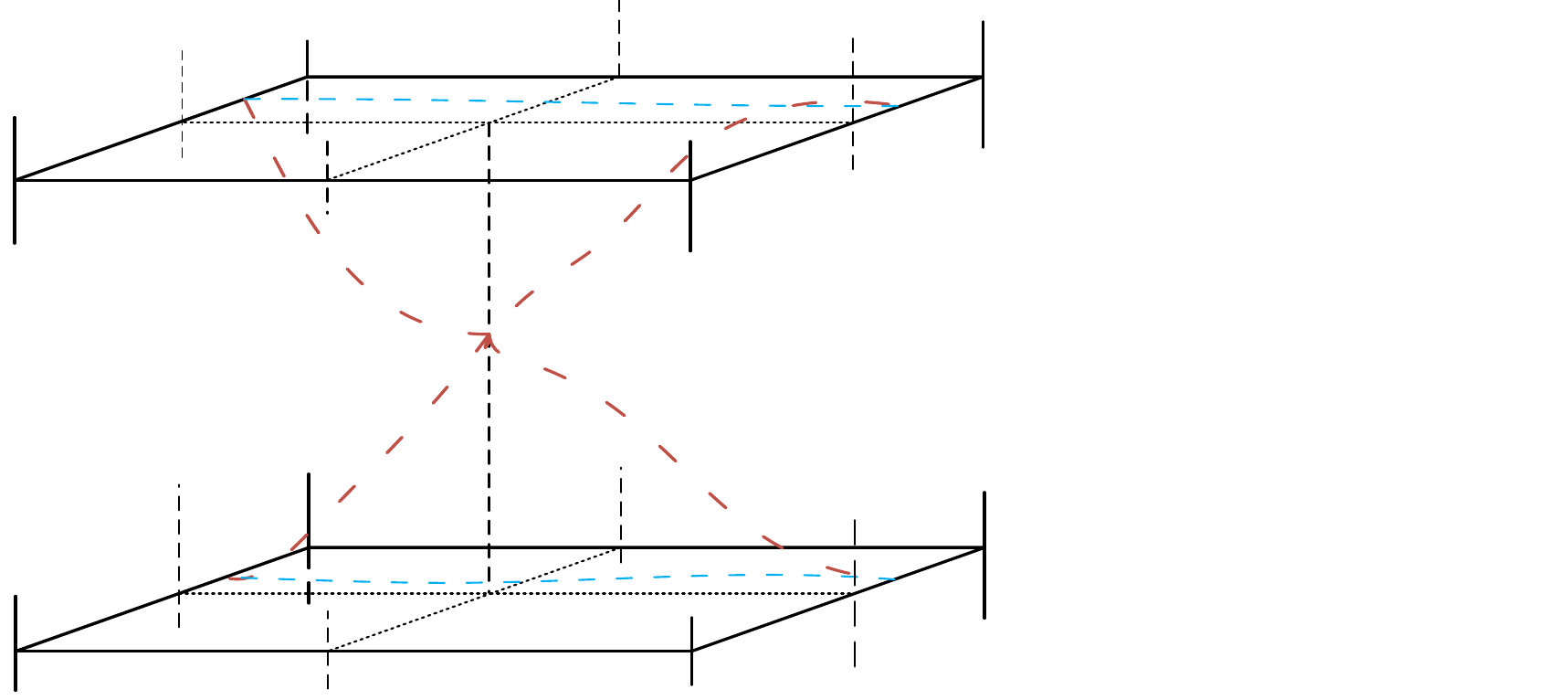}}%
    \put(0.39380441,0.26237756){\color[rgb]{0,0,0}\makebox(0,0)[lt]{\lineheight{1.25}\smash{\begin{tabular}[t]{l}$\alpha_{F_3,1}$\end{tabular}}}}%
    \put(0.42098079,0.17933649){\color[rgb]{0,0,0}\makebox(0,0)[lt]{\lineheight{1.25}\smash{\begin{tabular}[t]{l}$\alpha_{F_3, e_3}$\end{tabular}}}}%
    \put(0.1484175,0.17549113){\color[rgb]{0,0,0}\makebox(0,0)[lt]{\lineheight{1.25}\smash{\begin{tabular}[t]{l}$\alpha_{F_3,e_1e_3}$\end{tabular}}}}%
    \put(0.14303578,0.25542082){\color[rgb]{0,0,0}\makebox(0,0)[lt]{\lineheight{1.25}\smash{\begin{tabular}[t]{l}$\alpha_{F_3,e_1}$\end{tabular}}}}%
    \put(0.37976931,0.38264434){\color[rgb]{0,0,0}\makebox(0,0)[lt]{\lineheight{1.25}\smash{\begin{tabular}[t]{l}$\beta_{3,1}$\end{tabular}}}}%
    \put(0.59109705,0.35999838){\color[rgb]{0,0,0}\makebox(0,0)[lt]{\begin{minipage}{0.40695464\unitlength}\raggedright \end{minipage}}}%
  \end{picture}%
\endgroup%

\caption{closed circles in $M$}
\end{figure}
in {\bf Case 2}, we define
\begin{align*}
i_*:\pi_1(M_F,p_0)\longrightarrow& \pi_1(M,v)\\
 \beta_{i,g}\longmapsto& \alpha_{F_i,g}, \text{ for } i=1,2,4;\\
  \beta_{3,g}\longmapsto&  \alpha_{F_3,g\cdot \lambda(F_3)\lambda(F_1)}\alpha_{F_3,g}=\alpha_{F_3,ge_1e_3}\alpha_{F_3,g} 
\end{align*}
where $g\in\mathbb{Z}_2^3$. It can be check that the diagram (\ref{Diag1}) is communicate  both in  {\bf Case 1 }\&{ \bf 2}.

\begin{lemma}\label{lm1}
Both $i_*$ and $j_*$ are injective.
\end{lemma}
\noindent{\bf Proof.} The injectivivity of $j_*$ is showed in \cite[Theorem 3.3]{WY}, and the injectivivity of $i_*$ can be proved by the injectivivity of $j_*$ and the commutativity of diagram (\ref{Diag1}).  $\hfill{} \Box$

\newpage
\begin{theorem}\label{T3}
Let $M$ be a $3$-dimensional small cover over a simple polytope $P$, then $M$ is atoroidal if and only if there is no $4$-belt in $P$. 
\end{theorem}
\noindent{\bf Proof.} The Elliptization and Hyperbolization Theorems together imply that every atoroidal closed $3$-manifold is either spherical or hyperbolic, i.e. $\pi_1(M)$ is finite or infinite without subgroup $\cong \mathbb{Z}^2$. If there exists a $4$-belt in $P$, then according to Lemma \ref{lm1}, there is a $\pi_1$-injective torus or Klein bottle in $M$, which implies that there is a subgroup $\cong \mathbb{Z}^2$ in $\pi_1(M)$. This is a contradiction.

Conversely, we show that $M$ toroidal implies a $4$-belt in $P$. If $M$ toroidal, then there is a subgroup $\mathbb{Z}^2$ in $\pi_1(M)$. And $\pi_1(M)$ embeds in $W_P$, thus the subgroup $\mathbb{Z}^2$ embeds in $W_P$. In other words, there are two free commutable elements in $W_P$,  written as $x,y$. We assume that $x=s_1s_2\cdots s_m$ and $y=t_1t_2\cdots t_n$ are shortest expressions of $x$ and $y$, where $s_i,t_j$ are generators of $W_P$, and $m,n\geq 2$. Then $xy=yx$ implies $s_it_j=t_js_i, \forall i,j$, i.e. $F_{s_i}\cap F_{t_j}\neq\varnothing, \forall i,j$. 
Since $x$ is free, there exist $s, s'$ such that $F_{s}\cap F_{s'}=\varnothing$. And because $y$ is free, there also exist $t, t'$ such that $F_{t}\cap F_{t'}=\varnothing$ and $t,t',s,s'$ are four distinct generators. Otherwish, we have $(t,t')=(s,s')$ or $(t,t')=(s',s)$, now there must exist another two generators in $x$ or $y$ such that the intersection of the corresponding two faces is empty set.
Otherwish, the only two non-commutable gennerators in $x$ and $y$ is $s$ and $s'$, which implies that $x^2$ and $y^2$ is one of $(ss')^2$ and $(s's)^2$. Thus $x^2=y^2$ or $x^2=y^{-2}$, which contradict to $x,y$ generating $\mathbb{Z}^2$.
Hence there exists a $4$-belt.
$\hfill{} \Box$
 A simple polytope is {\em flag} if any pairwise intersecting faces have a common intersection. Davis \cite[Corollary 5.4]{D2} proved that a small cover is aspherical if and only if the orbit polytope is flag. And there are at most $3$ faces intersecting at one point in a $3$-dimensional simple polytope other than $\Delta^3$, So we have the proposition as following.
\begin{proposition}\label{P1}
Let $P\neq \Delta^3$ be a simple polytope of dimension $3$, then $P$ is flag if and only if there is no prismatic $3$-circuit in $P$.
\end{proposition}

\section{Andreev's Theorem in the right-angled case}

Andreev[1971] (see also \cite{RHD}) gives a complete characterization of compact hyperbolic polyhedra in dimension $3$ with nonobtuse dihedral angles. In the right-angled case, we have
\begin{theorem}[Andreev's Theorem in the right-angled case]\label{A-thm}
~\\A  simple polytope $P$($\neq \Delta^3$) has a geometric realization in $\mathbb{H}^3$ as a right-angled hyperbolic polytope if and only if there is no prismatic $3$ or $4$-circuit in $P$.
Furthermore, such geometric realization is unque up to isometry.
\end{theorem}

Now, we give a topological proof for the above theorem. For any $3$-dimensional simple polytope $P$, there exists a small cover $M$ over $P$ by the $4$-Colors Theorem. Conversely,
The group  $\mathbb{Z}_2^3$ that acts on a hyperbolic small cover of dimension $3$ produces
 a simple orbit polytope with all dihedral angles right-angled. Thus, we just need to prove that $M$ is hyperbolic if and only if there is no prismatic $3$ or $4$-circuit in $P$. \\
{\bf Proof of the necessity part.} If $M$ is hyperbolic, then $P$ can not be a $3$-simplex. If not, there is a element of order $2$ in $\pi_1(M)$, then $\pi_1(M)$ is not torsion-free, which is a contradiction. And other cases can be shown by Gauss-Bonnet theorem.
\\\\
{\bf Proof of the sufficiency part.}If there is no prismatic $3$-circuit, then $M$ is aspherical according to Proposition \ref{P1}. Furthermore, If there is no prismatic $3$-circuit, then there is no prismatic $4$-circuit if and only if  there is no $4$-belt in $P$.
 Thus $M$ is atoroidal by Theorem \ref{T3} if  there is no prismatic $3$ and $4$-circuit. Hence $M$ is hyperbolic by Thurston's Hyperbolization Theorem and such hyperbolic structure is unique up to isometry, thus $P$ has a unique geometric realization in $\mathbb{H}^3$, and now its all dihedral angles are right angle. $\hfill{} \Box$

\section*{Acknowledgments}
The author want to thank Yuxiu Lu, Hao Li and Jingfang Lian for their helpful comments.

\textsc{School of Mathematical Sciences, Fudan University, Shanghai, 200433, P.R.China.}

{\em E-mail address:} \url{wulisuwulisu@qq.com}
\end{document}